\documentclass[12pt, a4paper]{article}

\usepackage{pdfsync}
\usepackage{amsmath}\multlinegap=0pt
\usepackage[latin1]{inputenc}
\usepackage{amsthm}
\usepackage{amssymb,mathrsfs}
\usepackage{hyperref}
\usepackage{graphicx}
\usepackage{subcaption}
\usepackage{ragged2e}

\numberwithin{equation}{section}
\theoremstyle{plain}
\newtheorem{thm}{Theorem}[section]
\newtheorem{coro}[thm]{Corollary}
\newtheorem{prop}[thm]{Proposition}
\newtheorem{lem}[thm]{Lemma}

\newtheorem{assumption}[thm]{Assumption}
\theoremstyle{definition}

\theoremstyle{remark}
\newtheorem{rem}[thm]{Remark}

\newcommand{\M}{\mathcal{M}}
\newcommand{\R}{\mathbb{R}}
\newcommand{\Pro}{\mathbb{P}}

\newcommand\N{{\mathbb N}}

\newcommand\pref[1]{(\ref{#1})}

\let \eps\varepsilon

\DeclareMathOperator{\argmin}{argmin}

\DeclareMathOperator{\argmax}{argmax}
\DeclareMathOperator{\id}{id}

\newcommand\proj{\mathop{\mathrm{proj}}\nolimits}

\def\<#1,#2>{\left<#1,#2\right>}

\newcommand\ovv{{\overline v}}
\newcommand\ovh{{\overline h}}

\newcommand\undv{{\underline v}}
\newcommand\undh{{\underline h}}
\newcommand\ovnu{{\overline \nu}}

\def\PP{{\mathscr{P}}}

\title{On a Class of Adversarial Classification Problems which admit a Continuous Solution}

\author{Guillaume Carlier\thanks{CEREMADE, UMR CNRS 7534, Universit\'e Paris
Dauphine, PSL, Pl. de Lattre de Tassigny, 75775 Paris Cedex 16, FRANCE and INRIA-Paris, MOKAPLAN,
\texttt{carlier@ceremade.dauphine.fr}},
Maxime Sylvestre \thanks{CEREMADE, UMR CNRS 7534, Universit\'e Paris
Dauphine, PSL, Pl. de Lattre de Tassigny, 75775 Paris Cedex 16, FRANCE and INRIA-Paris, MOKAPLAN,
\texttt{maxime.sylvestre@dauphine.psl.eu}}
}

\begin{document}

\maketitle

\begin{abstract}  We consider a class of adversarial classification problems in the form of zero-sum games between a classifier and an adversary. The latter is able to corrupt data, at the expense of some optimal transport cost. We show that quite general assumptions on the classifier's loss functions and the adversary's transport cost functions ensure the existence of a Nash equilibrium with a continuous (or even Lipschitz) classifier's strategy. We also consider a softmax-like regularization of this problem and present numerical results for this regularization.\end{abstract}

\medskip

\textbf{MS Classification:} 49K35, 91A05, 49Q22.

\section{Introduction}

Despite their tremendous success   for performing difficult classification tasks (among others), machine learning algorithms raise major concerns about their robustness and vulnerability to small or even imperceptible perturbations of training data \cite{GoodfellowIntrig}. The problem of adversarial attacks and the analysis of perturbations of data leading to severe misclassification issues is now recognized as an important societal topic. This stimulated  in the last decade a huge stream of research which strongly suggests to replace the risk minimization paradigm which was prevalent in machine learning by alternative training methods ensuring robustness. Starting from the seminal paper \cite{NIPSGoodfellow2014},  adversarial training, i.e. training which is robust to (a class of) adversarial inputs, has gained a lot of popularity, see \cite{madry}, \cite{GoodfellowHarness} and the references therein.

\smallskip

Roughly speaking, the adversarial paradigm replaces risk minimization by its robust min-max (or worst case) counterpart where the max is taken over a set of adversarial scenarios. Such problems can be viewed as  zero-sum games between a classifier and an adversary.  The  game-theoretic viewpoint on adversarial training naturally leads to a number of questions such as the existence of a value, of optimal (possibly randomized, see \cite{Meunieretal}, \cite{pinot}) strategies, the characterization and computations of these strategies.... Of course, these issues depend very much on the learning problem at stake, the loss of the learner and the adversarial cost and admissible strategies. A reasonable framework for adversarial learning is the situation where the probability distribution of data can be slightly corrupted by the adversary. Understanding what a small perturbation of a distribution is or the cost of this perturbation, one enters the realm of \emph{optimal transport} theory. For instance, the recent series of works \cite{GT1}, \cite{GT2} and \cite{GTexist} demonstrated how robust multi-class classification can be reformulated in optimal transport terms. In \cite{GT2}, \cite{GTexist} as in \cite{pinot},   the adversary is allowed to move each individual data point up to certain  distance which makes his set of strategy a certain Wasserstein ball. The problem is also intimately related to the field of \emph{distributionally robust optimization}, \cite{Kuhn}, \cite{Blanchet}, where optimal transport has gained a lot of popularity as well, see \cite{Malick}.   Note also that optimal transport  distances can be used as  meaningful losses in the context of generative modeling as was shown in the influential work \cite{arjovsky} where Wasserstein GANs are introduced.

\smallskip

The first purpose of the present paper is to find, in the case of adversarial binary classification, simple conditions which guarantee not only that the game has a value but also that the classifier has a \emph{continuous} optimal strategy. Of course, in the case of hard transport constraints, such as in \cite{GT2}, \cite{pinot}, \cite{Meunieretal}, one cannot  hope that this is the case and even the existence of a Borel optimal classifier  is a delicate problem which was very nicely solved in \cite{GTexist}. This is why we will consider the softer situation where instead of a hard constraint the adversary has an integral optimal transport cost to pay to corrupt the data. We shall prove that continuity of this transport cost and mild convexity/monotonicity assumptions on the classifier's loss  ensure the existence of an optimal continuous strategy. A second contribution of our paper is the detailed investigation of the softmax regularization of the game and its use to solve some numerical examples. This approximation is of course  connected to entropic optimal transport theory which has gained increasing interest in the machine learning community thanks to Sinkhorn algorithm and the works of Cuturi and Peyr\'e, \cite{CuturiPeyre}, \cite{cuturilightspeed} \cite{cuturifast}.

\smallskip

The paper is organized as follows. The setting and assumptions are introduced in section \ref{sec-setting}, then both the classifier and the adversary's problems are reformulated in a concise way in section \ref{sec-concise}. In section \ref{sec-duality}, after showing that the game has a value, we show the existence of a continuous optimal strategy for the classifier. Section \ref{sec-softmax} is devoted to a detailed analysis of the softmax approximation of the classifier/adversary game and its convergence, several numerical examples based on this approximation  conclude the paper. Finally, section \ref{sec-concl} concludes  and gives some perspectives.

\section{Setting and Assumptions}\label{sec-setting} 

Let $(X, d)$ be a compact metric space equipped with its Borel $\sigma$-algebra, $(x,y)\in X\times \{-1, 1\}$ be random labeled points drawn according to some joint Borel probability distribution $\Pro$. The probability distribution $\Pro$ is  fully characterized by the two measures $A \subset X \mbox{ Borel} \mapsto \Pro(A\times \{i\})$, these are are not probability measures but coincide up to some normalization factors (that we will ignore in the sequel) with  the two probability distributions of $x$ given $y=i$ which we denote by $\mu_{i}\in \PP(X)$. 

\smallskip

The classifier has a loss function $(h,y)\in \R\times \{-1, 1\} \mapsto l(h,y)$ where $h$ is interpreted as a soft classifier for the binary label $y$. Defining $l_i:=l(\cdot, i)$ for $i=\pm 1$, we will always assume that $l_i$: $\R\to \R_+$ is convex (hence continuous) for $i=\pm 1$. Classical examples of such losses are
\begin{equation}\label{lossentr}
l(h,y)=\log(1+e^{-yh})
\end{equation}
and
\begin{equation}\label{lossmax}
l(h,y)=(1-yh)_+=\max(0, 1 -yh).
\end{equation}

The adversary may fool the classifier by perturbing the true  distributions $\mu_1$ and $\mu_{-1}$ at some cost. The adversary's (possibly randomized) strategies may be described by probability measures $\gamma_i \in \PP(X\times X)$ where $\gamma_i (A\times B)$ represent the probability that the true value of $x$ lies in $A$ but the observed (corrupted by the adversary) lies in $B$, given $y=i\in \{-1, 1\}$. The set of admissible strategies for the adversary thus reads\footnote{The notation $_\#$ will be used throughout the paper for pushforward of measures. More precisely, if $Z_1$ is a  metric space endowed with a Borel probability measure $m_1$, $Z_2$ is another metric space and $T$ is a measurable map $Z_1 \to Z_2$, $T_\#m_1$ is the probability measure  on $Z_2$ defined by $T_\# m_1(B)= m_1(T^{-1}(B))$  for every $B$, Borel subset of $Z_2$.}
\[\Delta_{\mu_1, \mu_{-1}} :=\{(\gamma_i)_{i=\pm 1} \in \PP(X\times X) \; : \; {\proj_1}_\# \gamma_i=\mu_i, \; i=\pm 1\}\]
where $\proj_1$ denotes the first projection\footnote{For $(x,z)\in X\times X$, $\proj_1(x,z)=x$ and  $\proj_2(x,z)=z$.} so that  ${\proj_1}_\# \gamma_i$ is the first marginal of $\gamma_i$. For each label $i$, there is a cost $c_i(x,z)$ incurred by the adversary for corrupting the true value $x$ into the corrupted value $z$, observed by the classifier. We  assume that  for  $i=\pm1$, the cost  $c_i$ : $X\times X \to \R_+ \cup\{+\infty\}$ is lsc,  and that  $c_i(x,x)=0$ for every $x\in X$. This implies that corrupting the data is  indeed costly for the adversary and the total cost of strategy $\gamma_1, \gamma_{-1}$ is given by
\[\sum_{i=\pm1} \int_{X\times X} c_i(x,z) \mbox{d} \gamma_i(x,z).\]

Typical examples of costs we have in mind are continuous costs of the form $c=d^\beta$ for some exponent $\beta>0$ or discontinuous (but lsc) costs like
\begin{equation}\label{costinfini}
c(x,z):=\begin{cases} 0 \mbox { if $d(x,z) \leq r$} \\+ \infty \mbox{ otherwise. }\end{cases}
\end{equation}

\begin{assumption}\label{ass:sec2}
To summarize our setting, throughout the paper we will always assume that
\begin{itemize}
\item $(X,d)$ is a compact metric space,

\item for $i=\pm 1$, the classifier's loss function $l_i$: $\R\to \R_+$ is convex,

\item  for $i=\pm 1$, the adversary's cost  $c_i$ : $X\times X \to \R_+ \cup\{+\infty\}$ is lsc and   $c_i(x,x)=0$ for every $x\in X$. 

\end{itemize}
\end{assumption}

All this results in the zero-sum game where the adversary's payoff is defined for every $h\in C(X)$ and $(\gamma_1, \gamma_{-1}) \in \Delta_{\mu_1, \mu_{-1}}$ by:
\begin{equation}\label{formedeJ}
 J(h, \gamma_1, \gamma_{-1}):= \sum_{i=\pm1} \int_{X\times X} [ l_i(h(z))-c_i(x,z) ] \mbox{d} \gamma_i(x,z).
 \end{equation}
Of course, setting
\[\ovv:=\inf_{h\in C(X)} \sup_{(\gamma_1, \gamma_{-1}) \in \Delta_{\mu_1, \mu_{-1}}}  J(h, \gamma_1, \gamma_{-1}) \]
and
\[\undv:=\sup_{(\gamma_1, \gamma_{-1}) \in \Delta_{\mu_1, \mu_{-1}}}   \inf_{h\in C(X)}  J(h, \gamma_1, \gamma_{-1}) \]
we have
\begin{equation}
\ovv \geq \undv.
\end{equation}
We now ask ourselves whether $\ovv=\undv$ i.e. the game has a value and whether there exist equilibrium strategies i.e. whether $J$ has a saddle-point $(h, \gamma_1, \gamma_{-1})$ in $C(X)\times \Delta_{\mu_1, \mu_{-1}}$. Of course, if $(h, \gamma_1, \gamma_{-1})$ is such a saddle-point then it is a Nash equilibrium in the sense that $h$ is an optimal classifier for the corrupted data $(\gamma_1, \gamma_{-1})$ and  $(\gamma_1, \gamma_{-1})$ maximizes the payoff of the adversary $J(h,\cdot,\cdot)$ over $ \Delta_{\mu_1, \mu_{-1}}$.

\section{Concise Formulations}\label{sec-concise}

The goal of this short section is to reformulate the inner maximization problem in the definition of $\ovv$ and the inner minimization problem in the definition of $\undv$.  Let us start by observing that for fixed $h\in C(X)$, and any $x\in X$, $i=\pm 1$, the set $\Gamma_i(x):=\argmax_{z\in X}\{ l_i (h(z))-c_i(x,z)\}$ is a nonempty compact subset of $X$ since $c_i$ is lsc and $l_i \circ h$ is continuous. Moreover, again by lower semicontinuity of $c_i$ and compactness of $X$, thanks to Proposition 7.33 in \cite{BertsekasShreve}, one can find a Borel map $S_i$ : $X\to X$ such that $S_i(x)\in \Gamma_i(x)$ for every $x\in X$. The maximization of $J(h,\cdot,\cdot)$ with respect to $\gamma_1, \gamma_{-1}$ in $\Delta_{\mu_1, \mu_{-1}}$ is therefore straightforwadly achieved by
\[\gamma_i:=(\id, S_i)_\# \mu_i.\]
We thus have 
\[\max_{(\gamma_1, \gamma_{-1}) \in \Delta_{\mu_1, \mu_{-1}}} J(h, \gamma_1, \gamma_{-1})=\sum_{i=\pm1} \int_X \max_{z\in X} \{l_i(h(z))-c_i(x,z)\} \mbox{d} \mu_i(x)\]
so that the upper value $\ovv$ reads as the value of the convex minimization problem:
\begin{equation}\label{upvalsimpl}
\ovv=\inf_{h\in C(X)} \sum_{i=\pm 1} \int_X \max_{z\in X} \{l_i ( h(z))-c_i(x,z)\} \mbox{d} \mu_i(x).
\end{equation}

As for the lower value, let us remark that it can be written in terms of the second marginals $\nu_i:={\proj_2}_\#\gamma_i$ as follows
\[\begin{split}
\undv= \sup_{(\gamma_1, \gamma_{-1}) \in \Delta_{\mu_1, \mu_{-1}}} \left\{ -\sum_{i=\pm 1} \int_{X\times X} c_i \mbox{d} \gamma_i + \inf_{h\in C(X)} \sum_{i=\pm 1} \int_X (l_i \circ h) \mbox{d} {\proj_2}_\# \gamma_i \right\}\\
= \sup_{(\nu_1, \nu_{-1}) \in \PP(X)^2}  \left\{ -\sum_{i=\pm 1} T_{c_i}(\mu_i, \nu_i) + \inf_{h\in C(X)} \sum_{i=\pm 1} \int_X (l_i \circ h) \mbox{d} \nu_i \right\}
\end{split}\]
where $T_{c_i}(\mu_i, \nu_i)$ denotes the value of the Monge-Kantorovich optimal transport problem:
\[T_{c_i}(\mu_i, \nu_i):=\min_{\gamma \in \Pi(\mu_i, \nu_i)}  \int_{X\times X} c_i(x,z) \mbox{d} \gamma(x,z)\]
where 
\[\Pi(\mu_i, \nu_i):= \{\gamma \in \PP(X\times X), \;  {\proj_1}_\# \gamma=\mu_i, \;  {\proj_2}_\# \gamma=\nu_i\}\]
denotes the set of transport plans between $\mu_i$ and $\nu_i$.  To further simplify the expression of $\undv$, it is useful to observe:

\begin{lem}\label{lemminenh}
Let $(\nu_1, \nu_{-1}) \in \PP(X)^2$ let $\ovnu:=\nu_1+\nu_{-1}$ and let $\alpha_i$ be the density of $\nu_i$ with respect to $\ovnu$ (so that $\alpha_{-1}=1-\alpha_1$ and $\alpha_1 \in [0, 1]$, $\ovnu$-almost everywhere), then
\begin{equation}\label{simpliinfh}
\inf_{h\in C(X)} \sum_{i=\pm 1} \int_X (l_i \circ h) \mathrm{d}  \nu_i= \int_X \Phi(\alpha_1(z)) \mathrm{d} \ovnu(z),
\end{equation}
where 
\begin{equation}\label{defdefi}
\Phi(\alpha):=\inf_{t\in \R} \{ \alpha l_1(t)+ (1-\alpha)l_{-1}(t)\}, \; \forall \alpha \in [0,1]. 
\end{equation}
\end{lem}

\begin{proof}
The fact that the left-hand side of \eqref{simpliinfh} is larger than its right-hand side follows from the very definition of $\Phi$. For the converse inequality, by monotone convergence, we first have 
\[\int_X \Phi(\alpha_1(z)) \mbox{d} \ovnu(z)=\lim_{M\to \infty } \int_X \Phi_M(\alpha_1(z)) \mbox{d} \ovnu(z),\]
where  
\[\Phi_M(\alpha):=\min_{t\in [-M,M]} \left\{ \alpha l_1(t)+ (1-\alpha)l_{-1}(t) + \frac{ t^2}{2M^4}\right\}, \; \forall \alpha \in [0,1].\]
Note that the minimum defining $\Phi_M(\alpha)$ is uniquely attained and that the point,  which we denote by $T^*_M(\alpha)$,  where this  minimum is achieved depends continuously on $\alpha$. Now remark that
\[ \int_X \Phi_M(\alpha_1(z)) \mbox{d} \ovnu(z) \geq  \sum_{i=\pm 1} \int_X l_i (T_M^*(\alpha_1(z))) \mbox{d} \nu_i(z)\]
and that, thanks to Lusin's Theorem (see \cite{Rudin}), for every $\eps>0$,  there exists $h_\eps\in C(X)$ such that $-M \leq h_\eps \leq M$ and $\ovnu(\{z\in X \; : \; h_\eps(z)\neq T_M^*(\alpha_1(z))\})\leq \eps$, so
\[ \begin{split}\int_X \Phi_M(\alpha_1(z)) \mbox{d} \ovnu(z) \geq  \sum_{i=\pm 1} \int_X l_i (h_\eps(z)) \mbox{d} \nu_i(z)-\eps \max_{[-M,M]}(l_1+l_{-1})\\
\geq  \inf_{h\in C(X)} \sum_{i=\pm 1} \int_X l_i \circ h \mbox{d} \nu_i-\eps \max_{[-M,M]}(l_1+l_{-1}). 
\end{split}\]
Letting $\eps\to 0$ and then $M\to \infty$, we  get the desired inequality and thus deduce \eqref{simpliinfh}.

\end{proof}

Thanks to lemma \ref{lemminenh}, the lower value can be rewritten as
\begin{equation}\label{simplelowerval}
\undv:=\sup_{(\nu_1, \nu_{-1}) \in \PP(X)^2}   -\sum_{i=\pm 1} T_{c_i}(\mu_i, \nu_i) +\int_X \Phi\Big( \frac{\mbox{d} \nu_1} {\mbox{d} ( \nu_1 +\nu_{-1})}  \Big) \mbox{d}(\nu_1+ \nu_{-1})
\end{equation}
 where $\Phi$ is the concave function defined in \eqref{defdefi}.

\section{Duality and Attainment Results}\label{sec-duality}

 \subsection{Duality}

\begin{prop}\label{minmax}
Under  assumption \ref{ass:sec2}, one has
\[ \ovv = \undv= \max_{(\gamma_1, \gamma_{-1}) \in \Delta_{\mu_1, \mu_{-1}}}   \inf_{h\in C(X)}  J(h, \gamma_1, \gamma_{-1}).\]
\end{prop}

\begin{proof}
We first claim that \eqref{upvalsimpl} can be reformulated as
\begin{equation}\label{reforavecfi}
\inf_ {(\varphi_1, \varphi_{-1}, h)\in C(X)^3} \left\{ \sum_{i=\pm 1} \int_X \varphi_i \mbox{d}\mu_i \; : \; \varphi_i(x)+c_i(x,z) \geq l_i(h(z)), \; \forall (x,z)\in X^2, \; i=\pm 1\right\}. 
\end{equation}
Note first that the inequality $\ovv \leq \inf \eqref{reforavecfi}$ is obvious. To show the converse inequality, take $h\in C(X)$ and set
\[\varphi_i(x):=\max_{z\in X} \{ l_i(h(z))-c_i(x,z)\}, \; x\in X, \; i=\pm 1.\]
Note now that $\varphi_i$ is usc and bounded. By standard arguments, $\varphi_i$ can be written as $\varphi_i=\inf_{n\in \N^*} \varphi_i^n$ where $\varphi_i^n$ is the $n$-Lipschitz function $x\in X\mapsto \varphi_i^n(x):=\max_{y\in X} \{\varphi_i(y)-n d(x,y)\}$, in particular, we have
\[\int_X \max_{z\in X} \{ l_i(h(z))-c_i(x,z)\} \mbox{d} \mu_i(x) =\lim_{n\to +\infty}  \int_X \varphi_i^n \mbox{d} \mu_i\]
since $(\varphi_1^n, \varphi_{-1}^n, h)$ is admissible for \eqref{reforavecfi} and $h\in C(X)$ is arbitrary, we can easily conclude that  $\ovv \geq \inf \eqref{reforavecfi}$. 

\smallskip

Let us now rewrite the program  \eqref{reforavecfi} in Fenchel-Rockafellar form as
\[\inf_ {(\varphi_1, \varphi_{-1}, h)\in C(X)^3} F(\varphi_1, \varphi_{-1}, h)+ G(\Lambda(\varphi_1, \varphi_{-1}, h))\]
where:
\begin{itemize}

\item $ F$ is the linear  form 
\[(\varphi_1, \varphi_{-1}, h)\mapsto F(\varphi_1, \varphi_{-1}, h):=  \sum_{i=\pm 1} \int_X \varphi_i \mbox{d}\mu_i,\] 

\item $\Lambda$ is the linear continuous operator from $C(X)\times C(X)\times C(X)$ to $C(X\times X)\times C(X\times X)\times C(X)$ defined through
\[(\varphi_1, \varphi_{-1}, h)\mapsto  \Lambda(\varphi_1, \varphi_{-1}, h)=(\varphi_1 \circ \proj_1, \; \varphi_{-1} \circ \proj_1, h),\]

\item and for $(\psi_1, \psi_{-1}, h)\in C(X\times X)^2 \times C(X)$,  
\[G(\psi_1, \psi_{-1},h)=\begin{cases} 0  \mbox{ if $\psi_i(x,z)+c_i(x,z)\geq l_i(h(z)), \; \forall (x,z)\in X^2, \; i=\pm1$,}  \\+\infty \mbox{ otherwise.} \end{cases}\]
\end{itemize}

Identifying the dual  of $C(X)$ (resp. $C(X\times X)$) with the space of measures on $X$  (respectively $X\times X$) which we denote by $\M(X)$ (resp. $\M(X\times X)$) , a  direct application of the Fenchel-Rockafellar duality Theorem (see \cite{EkelandTemam})  gives
\begin{equation}\label{rocka}
\ovv=\max_{(\gamma_1, \gamma_{-1}, \nu)\in \M(X\times X)^2 \times \M(X)} -F^*(\Lambda^*(\gamma_1, \gamma_{-1}, \nu))-G^*(-\gamma_1, -\gamma_{-1}, -\nu)
\end{equation}
where $\Lambda^*$ is the adjoint of $\Lambda$ which explicitly writes as 
\[\Lambda^*(\gamma_1, \gamma_{-1}, \nu)=({\proj_1}_\#\gamma_1, {\proj_1}_\#\gamma_{-1}, \nu), \; \forall (\gamma_1, \gamma_{-1}, \nu)\in \M(X\times X)^2 \times \M(X)\]
whence
\[F^*(\Lambda^*(\gamma_1, \gamma_{-1}, \nu))=\begin{cases} 0 \mbox{ if ${\proj_1}_\#\gamma_1=\mu_1, {\proj_1}_\#\gamma_{-1}=\mu_{-1}$, and $\nu=0$,}\\ +\infty \mbox{ otherwise } \end{cases}\]
and a direct computation gives 
\[-G^*(-\gamma_1, -\gamma_{-1}, 0)=\begin{cases} \inf_{h\in C(X)} \sum_{i=\pm 1} J(h, \gamma_1, \gamma_{-1})  \mbox{ if $\gamma_i \geq 0$, $i=\pm 1$} \\ -\infty \mbox{ otherwise.}\end{cases}\]
Replacing in \eqref{rocka} thus yields
\[\undv=\max_{(\gamma_1, \gamma_{-1})\in \Delta_{\mu_1, \mu_{-1}}} \inf_{h\in C(X)} J(h, \gamma_1, \gamma_{-1})=\ovv.\]
 
\begin{rem}
Note that the above application of the Fenchel-Rockafellar Theorem says that the $\sup$ in the $\sup \inf$ definition of $\undv$ is indeed a $\max$. The fact that this sup is achieved can  also be directly checked from the concise formulation in \eqref{simplelowerval} which is the maximization of a weakly $*$ usc function (for the term involving $\Phi$, this follows directly from Lemma \ref{lemminenh}) over the weakly $*$ compact set $\Delta_{\mu_1, \mu_{-1}}$.
\end{rem}

\end{proof}

\subsection{Existence of Continuous Optimal Classifiers}

It cannot be taken for granted for discontinuous $c_i$'s that the classifier has a continuous optimal classification strategy i.e. that  \eqref{upvalsimpl} admits solutions. However, under quite general conditions including the continuity of the costs $c_i$, we have:

\begin{thm}\label{existhc}
In addition to   assumption \ref{ass:sec2}, let us assume that
\begin{itemize}

\item $c_1$ and $c_{-1}$ are continuous on $X\times X$,

\item $l_1$ is nonincreasing, $l_{-1}$ is nondecreasing,

\item $\lim_{t\to -\infty}  l_1(t)=\lim_{t\to \infty}  l_{-1}(t)=+\infty$.
\end{itemize}
Then, \eqref{upvalsimpl} admits at least one solution, i.e. the classifier has a continuous optimal strategy. 

\end{thm}

\begin{proof}
Let $(h^n)_n\in C(X)^\N$ be a minimizing sequence for \eqref{upvalsimpl}, that is
\begin{equation}\label{minseq}
\lim_{n\to +\infty} \sum_{i=\pm 1} \int_X \varphi_i^n(x) \mbox{d} \mu_i(x)=\ovv
\end{equation}
where
\[\varphi_i^n(x):=\max_{z\in X} \{l_i(h^n(z))-c_i(x,z)\}, \; i=\pm 1, \, x\in X.\]
Since $c_1$  and $c_{-1}$ are uniformly continuous, the modulus
\[\omega(t):= \max_{i=\pm1} \max\{ \vert c_i(x,z)-c_i(x', z')\vert, \; (x,z,x',z')\in X^4, \; d(x,x')+d(z,z')\leq t\}\]
satisfies $\omega(t)\to 0$ as $t\to 0$. By definition of $\omega$, for every $(x,x',z)\in X^3$, we have
\[l_i(h^n(z))-c_i(x,z) \leq l_i(h^n(z))-c_i(x',z)+\omega(d(x,x'))\]
taking the supremum with respect to $z$ yields $\varphi_i^n(x) \leq \varphi_i^n(x') +\omega(d(x,x')$. Reversing the role of $x$ and $x'$, we obtain that 
 that for every $n$ and every $(x,x')\in X^2$ and $i=\pm1$, one has
\[\vert \varphi_i^n(x)-\varphi_i^n(x')\vert \leq \omega(d(x,x')).\]
So that both sequences $(\varphi_1^n)_n$ and $(\varphi_{-1}^n)_n$ are uniformly equicontinuous. Moreover since $l_i\geq 0$ and $c_i(x,x)=0$, $\varphi_i^n\geq 0$, but thanks to \eqref{minseq} we also have a uniform upper bound on $\min \varphi_i^n$, together with the uniform equicontinuity of $(\varphi_i^n)_n$ this shows that $(\varphi_i^n)$ is also uniformly bounded. Thanks to Arzel\`a-Ascoli Theorem, this shows that, up to a subsequence, $\varphi_i^n$ converges uniformly to some $\varphi_i$. Together with \eqref{minseq}, this entails
\begin{equation}\label{integralfiv}
\ovv=\sum_{i=\pm 1} \int_X \varphi_i(x) \mbox{d} \mu_i(x).
\end{equation}
Our goal now is to deduce the existence of an optimal solution to \eqref{upvalsimpl} which requires some more work. For $i=\pm 1$ and any $n$ and $z\in X$, set
\[\psi_i^n(z):=\min_{x\in X} \{\varphi_i^n(x) +c_i(x,z)\}\]
and observe that uniform convergence of $\varphi_i^n$ to $\varphi_i$ implies uniform convergence of $\psi_i^n$ to $\psi_i$ given by
\begin{equation}\label{psiconj}
\psi_i(z):=\min_{x\in X} \{\varphi_i(x) +c_i(x,z)\}.
\end{equation}

By construction, we have 
\begin{equation}\label{ineqgpsi}
\psi_1^n \geq l_1 \circ h^n, \; \psi_{-1}^n \geq l_{-1}\circ h^n,
\end{equation}
and since $\lim_{t\to -\infty}  l_1(t)=\lim_{t\to \infty}  l_{-1}(t)=+\infty$ and $\psi_i^n$ is uniformly bounded, this gives  a uniform bound on $h^n$: there is some $C$ such that $-C\leq h^n \leq C$ for every $n$. Now it is convenient to reexpress the inequalities in \eqref{ineqgpsi} in terms of $\psi_1^n$ and $\psi_{-1}^n$ only, using suitable generalized inverses of the convex and monotone losses $l_1$ and $l_{-1}$.  For $t\geq l_1(C)$, let
\[l_1^{-1}(t):=\inf\{ h\in (-\infty, C]  \; : \; l_1(h)\leq t\}\]
our assumptions ensure that $l_1^{-1}$ is nonincreasing continuous (recall that $l_1$ is both nonincreasing and convex so in fact, depending on whether its infimum is attained or not, it is either decreasing on $\R$, or decreasing on some interval  $(-\infty, a]$ and constant on $[a,+\infty)$) and that for $h\leq C$, the inequality $t\geq l_1(h)$ is equivalent to $h\geq l_1^{-1}(t)$. Similarly,   set for $t\geq l_{-1}(-C)$:
\[l_{-1}^{-1}(t):=\sup\{ h\in [-C, \infty)  \; : \; l_{-1}(h)\leq t\}\]
 so that $l_{-1}^{-1}$ is nondecreasing continuous and  for $h\geq- C$, the inequality $t\geq l_{-1}(h)$ is equivalent to $h\leq l_{-1}^{-1}(t)$. With these considerations in mind, \eqref{ineqgpsi} rewrites 
 \[l_{-1}^{-1}\circ \psi_{-1}^n \geq h^n \geq l_1^{-1} \circ \psi_1^n\]
 which, by continuity, passes to the limit as
\[\ovh:=l_{-1}^{-1} \circ \psi_{-1}  \geq \undh:= l_1^{-1} \circ \psi_1.\]
Note that both $\ovh$ and $\undh$ are continuous and that any $h\in C(X)$ such that $\undh\leq h\leq \ovh$ satisfies
\[\psi_i \geq l_i\circ h, \; i=\pm1.\]
Recalling \eqref{psiconj}, for such an $h$, we have
\[\max_{z\in X} \{l_i(h(z))-c_i(x,z)\} \leq \varphi_i(x), \; \forall x\in X\]
which, with \eqref{integralfiv}, implies that $h$ solves \eqref{upvalsimpl}.

\end{proof}

In fact, the above proof enables us to find even more regular solutions of the classifier problem:

\begin{prop}
If, in addition to the assumptions of Theorem \ref{existhc}, we assume that  
\begin{itemize}
\item either $c_1$  is Lipschitz on $X\times X$, and $l_1$ is either decreasing or of the form \eqref{lossmax},
\item or  $c_{-1}$  is Lipschitz on $X\times X$, and  $l_{-1}$ is either increasing or of the form \eqref{lossmax}.
\end{itemize}
Then, \pref{upvalsimpl} admits at least one Lipschitz solution.
\end{prop}

\begin{proof}
Following the proof of Theorem \ref{existhc}, we find $\psi_1$ and $\psi_{-1}$ of the form \eqref{psiconj} (that is $\psi_i$ is $c_i$-concave in the usual optimal transport terminology) such that  $l_1^{-1} \circ \psi_1 \leq l_{-1}^{-1} \circ \psi_{-1}$ and any $h\in C(X)$ such that $ l_1^{-1} \circ \psi_1 \leq h \leq l_1^{-1} \circ \psi_1$ solves \eqref{upvalsimpl}. Because $\psi_1$ and $\psi_{-1}$ are of the form \eqref{psiconj}, if $c_1$ (respectively $c_{-1}$) is Lipschitz so is $\psi_1$ (respectively $\psi_{-1}$). Assume that $c_1$ is Lipschitz. If $l_1$ is decreasing, since it is convex, its  subgradient is bounded away from zero on each compact interval, so that its inverse, $l_1^{-1}$ is in fact Lipschitz on compact sets and then so is $l_1^{-1}\circ \psi_1$. The  case where $l_1$ is of the form \eqref{lossmax} leads to the same conclusion since the generalized inverse $l_1^{-1}$ is also Lipschitz in this case. Finally, if  $c_{-1}$ is Lipschitz and $l_{-1}$ is increasing, or of the form \eqref{lossmax}, we similarly obtain a Lipschitz solution by  considering  $l_{-1}^{-1} \circ \psi_{-1}$ instead of $l_1^{-1}\circ \psi_1$.

\end{proof}

\begin{rem}.
If the losses attain their minimum and arrive there with a zero slope (for instance, with a quadratic behavior like $l_i(h)=(1-ih)_+^2$), their inverse functions are not Lipschitz (in the quadratic example $(1-ih)_+^2$, they behave like a square root near $0$ and are no better than $\frac{1}{2}$-H\"older continuous). 
\end{rem}
\begin{rem}
The previous proofs heavily rely on the continuity of the $c_i$'s, convexity and monotonicity of the $l_i$'s but the measures $\mu_i$ are totally arbitrary as well as the compact metric space $X$.
\end{rem}

\subsection{Primal-Dual Conditions and Consequences}

In this paragraph, we make the same assumptions as in Theorem \ref{existhc}. Under these assumptions, we know that both the classifier's problem \eqref{upvalsimpl} and the adversary's problem  \eqref{simplelowerval} admit solutions and that they have the same value by Proposition \ref{minmax}. We may take further advantage of duality  to obtain a characterization and some properties of solutions  by exploiting  primal-dual optimality conditions derived from Proposition \ref{minmax}. Let $h\in C(X)$, it directly follows from  Proposition \ref{minmax} that $h$ solves \eqref{upvalsimpl} if and only if there exist $(\nu_1, \nu_{-1})\in \PP(X)^2$ such that
\begin{align}
\sum_{i=\pm1} T_{c_i}(\mu_i, \nu_i)+\sum_{i=\pm 1} \int_X \max_{z\in X} \{l_i ( h(z))-c_i(x,z)\} \mbox{d} \mu_i(x) \nonumber\\
=\int_X \Phi\Big( \frac{\mbox{d} \nu_1} {\mbox{d} ( \nu_1 +\nu_{-1})}  \Big) \mbox{d} (\nu_1+ \nu_{-1}) \label{primaldual0}
\end{align} 
in which case $(\nu_1, \nu_{-1})$ automatically solves \eqref{simplelowerval}. Likewise $(\nu_1, \nu_{-1})\in \PP(X)^2$ solves  \eqref{upvalsimpl}  if and only if there exists $h\in C(X)$ such that  \eqref{primaldual0} holds, in which case $h$ automatically solves \eqref{upvalsimpl}. The basic primal-dual relation \eqref{primaldual0} therefore characterizes optimal strategies. We now recall (see \cite{Santambrogio}, \cite{Villani}) the Kantorovich duality formula which expresses the value of the optimal transport problem  $T_{c_i}(\mu_i, \nu_i)$ as 
\begin{equation}\label{dualTci}
 T_{c_i}(\mu_i, \nu_i)=\sup_{\psi \in C(X)}  \left\{ \int_X \psi^{c_i} \mbox{d} \mu_i + \int_X \psi \mbox{d} \nu_i  \right\}
\end{equation}
where $\psi^{c_i}$ denotes the $c_i$ transform of $\psi$:
\[\psi^{c_i}(x):= \min_{z\in X} \{c_i(x, z)-\psi(z)\}, \; \forall x\in X.\]
Optimizers in \eqref{dualTci} are called Kantorovich potentials between $\mu_i$ and $\nu_i$ and they are well-known to exist (if $c_i$ is continuous and $X$ is compact (see \cite{Santambrogio}, \cite{Villani}). We then have

\begin{prop}\label{primaldualprop}
Let $h\in C(X)$, $(\nu_1, \nu_{-1})\in \PP(X)^2$ set $\ovnu=\nu_1+ \nu_{-1}$ and let $\alpha_1$ be the density of $\nu_1$ with respect to $\ovnu$. Then $(h, \nu_1, \nu_{-1})$ satisfies the primal-dual relation \eqref{primaldual0} if and only if:
\begin{itemize}

\item for $i=\pm 1$, $l_i \circ h$ is a  Kantorovich potential between $\mu_i$ and $\nu_i$,

\item and, for $\ovnu$-a.e. $z\in X$, one has
\begin{equation}\label{primaldual2}
h(z)\in \argmin_{t\in \R}\{ \alpha_1(z) l_1(t)+ (1-\alpha_1(z)) l_{-1}(t)\}.
\end{equation}
\end{itemize}

\end{prop}

\begin{proof}

By Kantorovich duality formula \eqref{dualTci}, we have for $i=\pm 1$,
\begin{align} 
T_{c_i}(\mu_i, \nu_i)&+ \int_X \max_{z\in X} \{l_i ( h(z))-c_i(x,z)\} \mbox{d} \mu_i(x) =T_{c_i}(\mu_i, \nu_i) -\int_X (l_i\circ h)^{c_i} \mbox{d} \mu_i \nonumber\\
&\geq \int_X (l_i \circ h) \mbox{d} \nu_i =\int_X   (l_i \circ h) \frac{\mbox{d} \nu_i}  {\mbox{d} \ovnu} \;   \mbox{d} \ovnu \label{primaldual11}
\end{align}
summing these inequalities and using the definition of $\Phi$ (see \eqref{defdefi}) yields
\begin{equation}\label{primaldual12}
\sum_{i=\pm 1}T_{c_i}(\mu_i, \nu_i) -\sum_{i=\pm 1} \int_X (l_i\circ h)^{c_i} \mbox{d} \mu_i \geq \int_X \Phi( \alpha_1(z)) \mbox{d} \ovnu (z).
\end{equation}
Hence  \eqref{primaldual0}  is equivalent to having equality in \eqref{primaldual11} for $i=\pm 1$ (which means that $l_i \circ h$ is a  Kantorovich potential between $\mu_i$ and $\nu_i$) and having equality in \eqref{primaldual12} i.e. for  $\ovnu$-a.e. $z\in X$
\[ \Phi( \alpha_1(z))=\alpha_1(z) l_1(h(z)) + (1-\alpha_1(z)) l_{-1}(h(z))\]
which is exactly \eqref{primaldual2}.

\end{proof}

A first consequence is that as soon as the loss functions do not achieve their infimum, the adversary uses strategies that mix the two labeled classes $i=1$ and $i=-1$ in the sense that:

\begin{coro}\label{mixingcoro}
In addition to  the assumptions of Theorem \ref{existhc}, assume that $l_1$ is decreasing and $l_{-1}$ is increasing. Then, if $(\nu_1, \nu_{-1})$ solves  \eqref{simplelowerval} and we set $\ovnu=\nu_1+ \nu_{-1}$, we have
\[ \frac{\mathrm{d} \nu_1} {\mathrm{d} \ovnu } \in (0,1), \; \mbox{$\ovnu$-a.e.}.\]
\end{coro}

\begin{proof}
Again  let $\alpha_1$ be the density of $\nu_1$ with respect to $\ovnu$. 
We know from Theorem \ref{existhc} that there exists  $h\in C(X)$ solving  \eqref{upvalsimpl}. It then follows from \eqref{primaldual2} that for $\ovnu$-a.e. $z$ for which $\alpha_1(z)=0$, $h(z)$ minimizes $l_{-1}$ but the infimum of $l_{-1}$ is not achieved; this shows that, $\ovnu$-a.e., $\alpha_1>0$. In a similar way $\alpha_1$ has to be strictly less than $1$,  $\ovnu$-a.e. since otherwise $l_1$ would achieve its minimum.

\end{proof}

\begin{rem}
Actually the argument in the previous proof shows that if $l_1$ is decreasing  (respectively $l_{-1}$ is increasing) then the density of $\nu_1$ with respect to $\ovnu=\nu_1+ \nu_{-1}$ is 
$\ovnu$-a.e. strictly less than $1$ (respectively strictly positive). 

\end{rem}

Another consequence is the uniqueness of the optimal classifier's strategy on the support of the adversary's strategies:

\begin{coro}
In addition to  the assumptions of Theorem \ref{existhc}, assume that $l_1$ is decreasing, that $l_{-1}$ is increasing and that either $l_1$ or $l_{-1}$ is strictly convex. Let $(\nu_1, \nu_{-1})$ solve \eqref{simplelowerval} and $\ovnu=\nu_1+ \nu_{-1}$, then all solutions of \eqref{upvalsimpl} coincide on the support of $\ovnu$.
\end{coro}

\begin{proof}
We know from Corollary \ref{mixingcoro} that $\alpha_1$, the density of $\nu_1$ with respect to $\ovnu$, lies in $(0,1)$  $\ovnu$-a.e.; with \eqref{primaldual2} this implies that whenever $h$ solves \eqref{upvalsimpl}, for $\ovnu$-a.e. $z$, $h(z)$ is the minimizer of the strictly convex function $\alpha_1(z) l_1(\cdot)+ (1-\alpha_1(z)) l_{-1}(\cdot)$. Since $h$ is continuous, it is in fact uniquely determined on the support of $\ovnu$.
\end{proof}

\section{Softmax Regularization, Numerical Results}\label{sec-softmax}

\subsection{Softmax Approximation}

In order to provide an algorithm for the computation of the optimal classifier as well as the optimal adversarial strategy, we replace the maximum in the problem \eqref{upvalsimpl} by a softmax approximation. More precisely, we fix a reference measure $m\in \PP(X)$ \emph{with full support}, and for $\eps >0$,   $\psi \in C(X)$, $x\in X$ and $i\in \{-1,1\}$, we define:
\[\psi^{[\eps,c_i,m]}(x):= \eps \log \left(\int_X \exp\left(\frac{\psi(z)-c_i(x,z)}{\eps}\right) \mbox{d} m(z)\right) , \; x\in X\]
the softmax approximation of
\[\psi^{[0,c_i]}(x):= \max_{z\in X} \{\psi(z)-c_i(x,z)\}.\]
Since $m$ has full support and $c_i$ is continuous, it is a classical fact that $\psi^{[\eps,c_i,m]}$ converges uniformly (and obviously from below) to $\psi^{[0,c_i]}$ as $\eps \downarrow 0$. 
We  then consider the following regularized classifier's problem:
\begin{equation}\label{upvalsoftmaxreg}
     \ovv_\eps=\inf_{h\in C(X)}  T_\eps(h), \; T_\eps(h):= \sum_{i=\pm 1} \int_X (l_i \circ h)^{[\eps,c_i,m]} \mbox{d}\mu_i
    \end{equation}
 and recall that \eqref{upvalsimpl} reads
 \[\ovv=\inf_{h\in C(X)}  T_0(h), \; T_0(h):= \sum_{i=\pm 1} \int_X (l_i \circ h)^{[0,c_i]} \mbox{d}\mu_i.\]

       \begin{rem} 
 The convergence (in the sense of $\Gamma$-convergence\footnote{see the textbook \cite{Braides} for an overview of $\Gamma$-convergence.} of $T_\eps$ to $T_0$ in $C(X)$) of the regularized problem \eqref{upvalsoftmaxreg} to the original problem \eqref{upvalsimpl} is easy to see. Indeed, since we assumed that $m$ has full support, for each $h\in C(X)$, $T_\eps(h)$ converges to $T_0(h)$, which implies the $\Gamma-\limsup$ inequality: $\limsup_\eps T_\eps(h_\eps)\leq T_0(h)$ with the constant recovery sequence  $h_\eps=h$. As for the $\Gamma-\liminf$ inequality, take a sequence $(h_\eps)$ that converges to $h$ in $C(X)$. Note that for $\eps < \eps'$ we have $T_\eps(h_\eps) \geq T_{\eps'}(h_\eps)$ by the monotonicity of the softmax. Thus, since the right hand side converges to $T_{\eps'}(h)$ by taking the $\liminf$ we get $\liminf_\eps T_\eps(h_\eps) \geq T_{\eps'}(h)$. And we conclude that the $\Gamma-\liminf$ inequality $\liminf_\eps T_\eps(h_\eps)\geq T_0(h)$ holds by letting $\eps'$ tend to $0$. This in principle  ensures that the regularized problem is well-suited in order to compute an optimal classifier i.e. a minimizer of $T_0$. Indeed, $\Gamma$-convergence ensures convergence of values and minimizers, provided we can prove that $T_\eps$ admits minimizers $h_\eps$ (in $C(X)$) and that a family of minimizers $(h_\eps)_\eps$ remain in a compact set of $C(X)$. This is not obvious a priori and we will need several steps (relaxation, duality and uniform continuity estimates) to achieve this goal.
 \end{rem}   
  
{\textbf{Relaxation.}} The existence of a minimizer of $T_\eps$ in $C(X)$  is not completely straightforward at first glance. We shall therefore  relax \eqref{upvalsoftmaxreg} to the larger class of $m$-measurable functions. Since $l_i\geq 0$ and $c_i$ is bounded,  $(l_i \circ h)^{[\eps,c_i,m]}$ is bounded from below as soon as $h$ is $m$-measurable, and a necessary and sufficient condition to prevent that it takes the value $+\infty$ somewhere (equivalently everywhere) on $X$ is the integrability of $\exp\Big({\frac{l_i \circ h}{\eps}}\Big)$. In other words, one can define $(l_i \circ h)^{[\eps,c_i,m]}$ for $h$ measurable and  in this  larger class of classifiers, the set where $T_\eps$ is finite is the (convex) set
\begin{equation}\label{defdeFeps}
F_\eps:=\left\{h \in L^{1}(m) \; :  \; \sum_{i=\pm1} \exp \Big({\frac{l_i\circ h}{\eps}}\Big) \in L^1(m)\right\}
\end{equation}
for $h$ in this class, not only $(l_i \circ h)^{[\eps,c_i,m]}$ is well defined, but it is actually continuous with a modulus of continuity which only depends on that of $c_i$ (see \cite{gerodim}). The relaxed formulation of \eqref{upvalsoftmaxreg} then reads: 
\begin{equation}\label{abasmacron}
\inf_{h\in F_\eps}  T_\eps(h).
\end{equation}  
 We then have
 
 \begin{prop}\label{existhlreg}
     Under the assumptions of Theorem \ref{existhc}, one has relaxation and existence of a relaxed solution to \eqref{abasmacron} in the sense that
     \[  \ovv_\eps=\min_{h\in F_\eps}  T_\eps(h).\]
  If one further assumes that either $l_1$ or $l_{-1}$ is strictly convex then the relaxed problem \eqref{abasmacron} admits  a unique solution.
\end{prop}
\begin{proof}
Of course, $\ovv_\eps$ is larger than or equal to the value of \eqref{abasmacron}. To prove the converse inequality, let  $h\in F_\eps$ and consider its truncations
 \[h_n:=\max(-n, \min(h, n)), \;  n\in \N^*.\]
  By monotonicity of $l_1$, $l_{-1}$, one has 
  \[\exp \Big({\frac{l_i\circ h_n}{\eps}} \Big)  \leq \max\Big( \exp \Big({\frac{l_i(0)}{\eps}}\Big), \exp\Big({\frac{l_i\circ h}{\eps}}\Big) \Big)\in L^1(m),\]
   so that one easily deduces from Lebesgue's dominated convergence Theorem that $T_\eps(h_n)$ converges to $T_\eps(h)$. Since each $h_n$ belongs to $L^\infty(m) \subset F_\eps$ this implies that the value of \eqref{abasmacron} coincides with $\inf_{L^\infty(m)}  T_\eps$. Now, for $h\in L^{\infty}(m)$ and $\delta >0$, by Lusin's Theorem, we can find $h_\delta \in C(X)$ with $\Vert  h_\delta \Vert_{L^\infty(m)} =  \Vert  h \Vert_{L^\infty(m)}$ such that $X_\delta:=\{h=h_\delta\}$ satisfies $m(X\setminus X_{\delta})\leq \delta$. For fixed $x$ and $i$, we then have for some positive constant $C$ (depending on $\Vert  h \Vert_{L^\infty(m)}$, $\eps$, $l_i$ and $c_i$ but not on $\delta$)
\[
\int_{X} \exp \Big({\frac{l_i(h_{\delta}(z))-c_i(x,z)}{\eps}} \Big) \mbox{d} m(z) \leq \int_{X} \exp \Big({\frac{l_i( h(z))-c_i(x,z)}{\eps}} \Big) \mbox{d} m(z)+ C\delta\]
from which we deduce that $T_\eps(h) \geq \limsup_{\delta \to 0} T_\eps(h_\delta) \geq \ovv_\eps$. Since $h\in L^\infty(m) $ is arbitrary in the previous argument, we obtain
\[\ovv_\eps \leq \inf _{L^\infty(m)} T_\eps=\inf_{F_\eps} T_\eps.\]

\smallskip

Let us now prove that $T_\eps$ admits a minimizer over $F_\eps$ (equivalently over measurable classifiers). One easily deduces from Fatou's Lemma (together with the nonnegativity of $l_i$ and the boundedness of $c_i$) that $T_\eps$ is lsc for the strong topology of $L^1$, since it is convex it is also weakly sequentially lsc for the weak topology of $L^1$.  Let $(h_n)$ be a minimizing sequence for \eqref{abasmacron}. The sequence of (nonnegative)  functions $(l_i \circ h_n)^{[\eps,c_i,m]}$ is bounded in $L^1(\mu_i)$, since it is also uniformly equicontinuous, it is uniformly bounded hence relatively compact in $C(X)$ by Arzel\`a-Ascoli's Theorem. But having a uniform bound on $(l_i \circ h_n)^{[\eps,c_i,m]}$ and using again the boundness of $c_i$  gives a uniform bound on $\int_X \sum_{i=\pm 1} e^{\frac{l_i \circ h_n}{\eps}} \mbox{d} m$. Now  observe that since $l_{-1}$ is convex nondecreasing and $\lim_{+ \infty } l_{-1}=+\infty$ one can find $a \in (0,1)$  such that $l_{-1}(t)\geq a \max(t,0)- \frac{1}{a}$ for every $t\in \R$. In a similar way, since $l_{1}$ is convex nonincreasing and $\lim_{- \infty } l_{1}=+\infty$ and up to choosing a possibly smaller $a$, we may also assume  $l_1(t)\geq a \max(-t,0)- \frac{1}{a}$ for every $t\in \R$. This implies that 
\[\max(l_1(t), l_{-1}(t)) \geq a \vert t \vert -\frac{1}{a}, \; \forall t\in \R,\]
which  gives a bound on $\int_X \exp\Big( \frac{a \vert h_n \vert}{\eps} \Big) \mbox{d} m$. In particular $(h_n)$ is uniformly integrable and therefore by the Dunford-Pettis Theorem, it has a subsequence that converges weakly in $L^1$ to some $h$. Thanks to the sequential weak lower semicontinuity in $L^1$ of $T_\eps$, we have $h\in F_\eps$ and $T_\eps(h) \leq \ovv_\eps$ so that $h$ solves \eqref{abasmacron}.  Finally, the uniqueness of a solution to \eqref{abasmacron} is a direct consequence of the strict convexity of $T_\eps$ in the case of the strict convexity of either $l_1$ or $l_{-1}$.
\end{proof}

{\textbf{Duality.}}  Integral functionals involving softmax as in the definition of $T_\eps$ are well-known to be related by convex duality to entropic optimal transport (see \cite{Nutz} and \cite{Leonard}). The regularized optimal transport cost between $\mu_i$ and $\nu_i$ with respect to the cost $c_i$ and the regularization parameter $\eps$ and the reference measure $m$ is defined by
\begin{equation}
     T_{c_i}^{\eps,m}(\mu_i, \nu_i) = \inf_{\gamma_i \in \Pi(\mu_i,\nu_i)} \int c_i d\gamma_i  + \eps H(\gamma_i | \mu_i \otimes m).
\end{equation}
where 
\[ H(\gamma_i | \mu_i \otimes m)=\begin{cases} \int_{X\times X} \log\Big( \frac{\mathrm{d} \gamma }{ \mathrm{d} (\mu_i \otimes m)} \Big) \mbox{d} \gamma \mbox{ if $\gamma\ll \mu_i \otimes m$} \\ +\infty \mbox{ otherwise } \end{cases}\]
denotes the relative entropy of $\gamma$ with respect to $\mu_i \otimes m$. This terms forces the (unique by strict convexity) optimal entropic plan to have a density with respect to $\mu_i \otimes m$ so that $T_{c_i}^{\eps,m}(\mu_i, \nu_i)$ is finite only if the second marginal $\nu_i$ is itself absolutely continuous with respect to $m$. The dual formulation of the  previous entropic optimal transport problem takes the form
\begin{equation}\label{dualeot}
 T_{c_i}^{\eps,m}(\mu_i, \nu_i) =\sup _{\psi \in C(X)} \left\{ \int_X \psi \mbox{d} \nu_i -\int_X  \psi^{[\eps, c_i, m]} \mbox{d} \mu_i \right\}
\end{equation}
If $\psi$ is a solution\footnote{In \eqref{dualeot}, one can consider more general potentials  than continuous ones, and maximize in the class of $m$-measurable $\psi$'s for which $\psi^{[\eps, c_i, m]} $ is finite i.e. satisfy the exponetial integrability condition $\exp(\eps^{-1} \psi) \in L^1(m)$.} of this dual formulation, it is called a \emph{Schr\"odinger potential} between $\mu_i$ and $\nu_i$.  Thanks to Fenchel-Rockafellar duality, arguing as in the proof of Proposition \ref{minmax}, the value $\ovv_\eps$ of  the regularized classifier problem \eqref{upvalsoftmaxreg} can be expressed as 
\begin{equation}\label{lowvalsoftmaxreg}
     \undv_\eps =\sup_{(\nu_1, \nu_{-1}) \in \PP(X)^2}   -\sum_{i=\pm 1} T_{c_i}^{\eps,m}(\mu_i, \nu_i) +\int_X \Psi\Big( \frac{\mbox{d} \nu_1} {\mbox{d} m}(z),\frac{\mbox{d} \nu_{-1}}{\mbox{d} m}(z)  \Big) \mbox{d}m(z)
\end{equation}
where $\Psi$ is defined similarly to $\Phi$ but here the reference measure is $m$ instead of $\nu_1+\nu_{-1}$ and thus the function depends on two variables $(\alpha_1,\alpha_{-1})$ which represent respectively the density of $\nu_1$ and $\nu_{-1}$ with respect to $m$. The precise definition of $\Psi$ is
\begin{equation*}
     \Psi(\alpha_1, \alpha_{-1}) = \min_{t\in \R}  \{\alpha_1 l_1(t) + \alpha_{-1} l_{-1}(t)\}.
\end{equation*}
Note that the existence of a solution to \eqref{lowvalsoftmaxreg} is guaranteed by the Fenchel-Rockafellar Theorem and that the equality $\ovv_\eps=\undv_\eps$ (absence of duality gap) can be interpreted as the existence of a value for a zero-sum game between a classifier and an adversary incurring an entropic optimal transport cost. 

\smallskip

As in the non regularized case, we can derive properties of the optimizers by exploiting the primal-dual optimality conditions. In particular, these optimality conditions allow us to show that the optimal classifier is continuous and thus that the original problem admits a solution. Let $h \in F_\eps$ (recall that $F_\eps$, defined in \eqref{defdeFeps} is the subset of $L^1(m)$ where $T_\eps$ is finite)  it directly follows from the absence of duality gap that $h$ solves \eqref{abasmacron} if and only if there exist $(\nu_1, \nu_{-1})\in \PP(X)^2$ such that
\begin{equation}\label{primaldualreg1}
     \sum_{i=\pm 1} T_{c_i}^{\eps,m}(\mu_i, \nu_i) + \int_X (l_i \circ h)^{[\eps,c_i,m]} \mbox{d} \mu_i(x)
     = \int_X \Psi\Big( \frac{\mbox{d} \nu_1} {\mbox{d} m}(z) ,\frac{\mbox{d} \nu_{-1}} {\mbox{d} m}(z)  \Big) \mbox{d}m(z)
\end{equation}
and likewise for $(\nu_1, \nu_{-1})$ which would solve \eqref{lowvalsoftmaxreg} if such an $h$ exists. We then have the following proposition.
\begin{prop}\label{primaldualregprop}
Let $h \in F_\eps$, $(\nu_1, \nu_{-1})\in \PP(X)^2$  and let $\alpha_i$ be the density of $\nu_i$ with respect to $m$. Then $(h, \nu_1, \nu_{-1})$ satisfies the primal-dual relation \eqref{primaldualreg1} if and only if:
\begin{itemize}
     \item for $i=\pm 1$, $l_i \circ h$ is a  Schr\"odinger potential between $\mu_i$ and $\nu_i$,
     \item for $m$-a.e. $z\in X$, one has
     \begin{equation}\label{primaldualreg11}
          h(z)\in \argmin_{t\in \R}\{ \alpha_1(z) l_1(t)+ \alpha_{-1}(z) l_{-1}(t)\},
     \end{equation}
     \item for $i=\pm 1$, $\nu_i$ satisfies,
     \begin{equation}
          \frac{\mathrm{d} \nu_i}{\mathrm{d} m}(z) = \int_X \frac{\exp\left(\frac{l_i ( h(z))-c_i(x,z)}{\eps}\right)}{\int_X \exp\left(\frac{l_i ( h(z'))-c_i(x,z')}{\eps}\right) \mathrm{d}m(z')} \mathrm{d}\mu_i(x).
     \end{equation}
\end{itemize}
\end{prop}
\begin{proof}
     The proof is similar to the one of proposition \ref{primaldualprop}. Indeed given $h\in F_\eps$, thanks to \eqref{dualeot} we have the following inequality:
     \begin{equation}
          T_{c_i}^{\eps,m}(\mu_i, \nu_i) + \int_X (l_i \circ h)^{[\eps,c_i,m]}  \mbox{d} \mu_i(x) \geq \int_X l_i\circ h \mbox{d} \nu_i
     \end{equation}
     which by summation gives
     \begin{equation}
          \sum_{i = \pm 1} T_{c_i}^{\eps,m}(\mu_i, \nu_i) + \int_X (l_i \circ h)^{[\eps,c_i,m]}  \mbox{d} \mu_i(x) \geq \int_X \Psi\Big( \frac{\mbox{d} \nu_1} {\mbox{d} m}(z),\frac{\mbox{d} \nu_{-1}} {\mbox{d} m}(z)  \Big)  dm(z).
     \end{equation}
     Note that \eqref{primaldualreg1} implies equality in the chain of inequalities which proves the desired properties of $h$. Now \eqref{dualeot} grants that the entropic optimal transport plan between $\mu_i$ and $\nu_i$ has a density with respect to $\mu_i \otimes m$ of the form
     \begin{equation}
          \frac{\mathrm{d}\gamma}{\mathrm{d} \mu_i \otimes m}(x,z) = \frac{\exp\left(\frac{l_i ( h(z))-c_i(x,z)}{\eps}\right)}{\int_X \exp\left(\frac{l_i ( h(z'))-c_i(x,z')}{\eps}\right) \mathrm{d} m(z')}.
     \end{equation}
     Integrating over $x$ with respect to $\mu_i$ gives back the second marginal of $\gamma$ and thus $\nu_i$ has the following density with respect to $m$:
     \begin{equation}
          \frac{\mathrm{d}\nu_i}{\mathrm{d} m}(z) = \int_X \frac{\exp\left(\frac{l_i ( h(z))-c_i(x,z)}{\eps}\right)}{\int_X \exp\left(\frac{l_i ( h(z'))-c_i(x,z')}{\eps}\right) \mathrm{d} m(z')} \mathrm{d}\mu_i(x).
     \end{equation}
\end{proof}

From the previous primal-dual optimality conditions, we can prove that the optimal classifier for the relaxed problem \eqref{abasmacron} is in fact continuous.
\begin{coro}\label{existhcreg}
    In addition to the assumptions of Theorem \ref{existhc}, assume that that $m$ has full support and either $l_1$ is decreasing or $l_{-1}$ is increasing. If $h \in F_\eps$ solves the relaxed problem \eqref{abasmacron}, it has a continuous representative. Hence the initial problem \eqref{upvalsoftmaxreg} admits a (unique) solution.
\end{coro}
\begin{proof}
     By Proposition \ref{primaldualregprop}, we know that, defining  $\alpha_i (z) = \frac{d\nu_i}{dm}(z)$, for $m$-a.e. $z\in X$, $h(z)$ is a minimizer of the function
     \begin{equation}
          g: \; s \mapsto \alpha_1(z)l_1(s)+ \alpha_{-1}(z)l_{-1}(s).
     \end{equation}
     This function is convex and finite hence right and left differentiable everywhere and thus by the characterization of minimizers we have 
     \begin{equation}\label{optigh}
     g'(h(z)^-) \leq 0 \leq g'(h(z)^+).
     \end{equation}
      Note that, again by Proposition \ref{primaldualregprop}, $\alpha_i(z) = \exp(\frac{l_i(h(z))}{\eps}) F_i(z)$ where 
        \begin{equation}\label{defFii}
          F_i(z) = \int_X \frac{\exp\left(\frac{-c_i(x,z)}{\eps}\right)}{\int_X \exp\left(\frac{l_i ( h(z'))-c_i(x,z')}{\eps}\right)\mathrm{d} m(z')} \mathrm{d}\mu_i(x),
     \end{equation}     
so that $F_i$ is a positive and continuous function. Let us now observe that
 the optimality conditions \eqref{optigh} satisfied by $h(z)$ are the same as for the minimization with respect to $s$  of the function
     \begin{equation}
          H : (s,z) \mapsto \sum_{i = \pm 1} \exp \Big(\frac{ l_i(s)}{\eps} \Big) F_i(z) 
     \end{equation}
     which, under our assumptions, is strictly convex and  coercive in $s$ (uniformly in $z$) and jointly continuous in $(s,z)$. This guarantees that  the map 
   $\tilde{h}: \; z \in X \mapsto \argmin_{s \in \mathbb{R}} H(s,z)$ is continuous. Now $h$ coincides $m$-a.e. with $\tilde{h}$ and thus the continuous function $\tilde{h}$ is also optimal and is the unique minimizer of \eqref{upvalsoftmaxreg} because $m$ has full support. 
   
\end{proof}
Finally in order to have a quantitative rate of convergence of the value of the regularized problem to the value of the original problem we strengthen the assumptions and show that the optimal classifier is Lipschitz uniformly in the regularization parameter. 
\begin{prop}\label{hreglip}
 In addition to the assumptions of Theorem \ref{existhc}, assume that $m$ has full support and that
      \begin{itemize}
          \item for $i = \pm 1$, $l_i$ is differentiable,  $l_1$ is decreasing and  $l_{-1}$ is increasing,
          \item for $i = \pm 1$,  $c_i$ is Lipschitz in the second variable uniformly in the first one.
     \end{itemize}
     Then, the optimal classifier i.e. the solution of \eqref{upvalsoftmaxreg} is bounded and Lipschitz uniformly in $\eps$.
\end{prop}
\begin{proof}

     Let $\eps>0$ and $h$ be the solution of \eqref{upvalsoftmaxreg}. In what follows $C$ will denote a positive constant that does not depend on $\eps$ but which may change from a line to another. Using the same notations as in  Corollary \ref{existhcreg}  and its proof, we know that $h(z)$ is the  minimizer of 
     \begin{equation}
           g_z: \; s \mapsto \sum_{i = \pm 1} \exp\left(\frac{l_i(s)}{\eps}\right)\exp\left(\frac{K_i(z)}{\eps}\right)
     \end{equation}
where $K_i(z) = \eps \log\left(F_i(z) \right)$ (with $F_i$ defined by \eqref{defFii}). That is $h(z)$ is the unique root of the equation  
     \begin{equation*}
         \sum_{i=\pm 1} l_i'(s)\exp\left(\frac{l_i(s)}{\eps}\right)\exp\left(\frac{K_i(z)}{\eps}\right) = 0.
     \end{equation*}
            Remark  that by boundedness of $c_i$ and optimality of $h$ we have 
     \[ \sum_{i=\pm 1} \eps \log \Big(\int_X \exp \Big(\frac{l_i\circ h}{\eps} \Big) \mathrm{d} m \Big) \leq C\]
     this implies that $K_i(z)$  is bounded uniformly in $z$ and in $\eps$. Together with our Lipschitz assumption on $c_i$, this gives
     \begin{equation}\label{uboundK}
      \vert K_i (z)\vert \leq C, \; \vert K_i(z_1)-K_i(z_0) \vert \leq C \vert z_1-z_0\vert, \; \forall (z, z_0, z_1) \in X^3.
      \end{equation}
     Moreover since $s=h(z)$ minimizes $g_z$, we have 
     \begin{equation}
          \sum_{i = \pm 1} \exp\left(\frac{l_i(s)}{\eps}\right)\exp\left(\frac{K_i(z)}{\eps}\right) \leq \sum_{i =\pm 1}\exp\left(\frac{l_i(0)}{\eps}\right)\exp\left(\frac{K_i(z)}{\eps}\right).
     \end{equation}
     Taking the logarithm of each side of the previous inequality, by  \eqref{uboundK} and concavity of the logarithm, we obtain
     \begin{equation}\label{boundsk}
          \sum_{i = \pm 1} l_i(s) + K_i(z) \leq 2 \max_{i=\pm1} l_i(0) + 2 \max_{i=\pm 1} K_i(z)\leq C.
     \end{equation}
     We already observed in the proof of Proposition \ref{existhlreg} that $l_1(s)+l_{-1}(s) \geq a\vert s\vert-\frac{1}{a}$ for some $a\in (0,1)$ and every $s$, hence $h$ can be uniformly bounded independently of $\eps$.
     
       \smallskip
    
    Let us now show that $h$ is Lipschitz on $X$ uniformly in $\eps$. Let us further assume for the moment that $l_1$ and $l_{-1}$ are of class $C^2$.  Take $z_0,z_1 \in X$ define $s_0=h(z_0)$, $s_1=h(z_1)$ and for $t\in [0,1]$, set $K_i^t:=(1-t)K_i(z_0)+ tK_i(z_1)$ and
    \[G_t(\cdot):= \sum_{i=\pm 1} \exp\left(\frac{l_i(\cdot)+K_i^t}{\eps} \right), \; s(t) :=\argmin G_t\]
    so that $s(t)$ is a curve joining $s_0$ to $s_1$. Arguing as for \eqref{boundsk} and using the uniform bound \eqref{uboundK}, we can find a uniform bound for $s(t)$: $s(t)\in [-M, M]$ for every $t\in [0,1]$. It follows from the implicit function Theorem that $s(\cdot)$ is differentiable and by differentiating the optimality condition $G'_t(s(t))=0$, we see that its derivative $\dot{s}(t)$ satisfies
    \[ \vert \dot{s}(t) \vert \leq \frac{ \sum_{i=\pm 1} e^{\frac{l_i(s(t))+K_i^t}{\eps}} \vert l'_i(s(t))\vert \vert K_i(z_1)-K_i(z_0)\vert  } { \sum_{i=\pm 1} e^{\frac{l_i(s(t))+K_i^t}{\eps}} [\eps l''_i(s(t))+ l'_i(s(t))^2]}\]
    hence by \eqref{uboundK} and using the fact that $\vert l'_i\vert$ is bounded and bounded away from $0$ on $[-M,M]$, we get
    \[\vert h(z_1)-h(z_0)\vert\ \leq \int_0^1 \vert  \dot{s} \vert \leq C \max_{i=\pm 1}  \vert K_i(z_1)-K_i(z_0) \vert \leq C \vert z_1-z_0\vert\]
    where the constant $C$ depends on the first derivatives of the losses but neither on their second derivatives nor on $\eps$. If the losses $l_i$ are just differentiable, we can approximate them by $C^2$ ones by convolution. By the previous argument, the optimal classifiers for these approximated losses are uniformly Lipschitz (with a Lipschitz constant that does not depend on $\eps$) and converge  (by uniqueness) to the optimal classifier for the initial losses $l_i$.  We can therefore remove the extra assumption that $l_i$ is $C^2$ and reach the very same estimate which enables us to conclude that the solution of \eqref{upvalsoftmaxreg} is Lipschitz uniformly in $\eps$. 
    \end{proof}
We conclude with the following rate of convergence which is a direct consequence of the rate of convergence of the softmax towards the maximum. 
\begin{prop}
Under the assumptions of Proposition \ref{hreglip}, if we further assume that there is $\delta > 0$ such that 
\begin{equation}\label{ahlfors}
     \liminf_{r \to 0} \inf_{x\in Z}\frac{m(B(x,r))}{r^{\delta}} > 0.
\end{equation}
Then we have $0 \leq \ovv -\ovv_\eps  =O( \eps \log(\frac{1}{\eps}))$ as $\eps \downarrow 0$.
\end{prop}
\begin{proof}
The lower bound is a direct consequence of the softmax being smaller than the maximum. For the upper bound, we will use the following argument. Given $L>0$ and  an $L$-Lipschitz function $g$  from $X$ to $\R$, we have for $\eps >0$:
\begin{equation*}
     \eps \log\left(\int_X e^{g(z)/\eps} \mbox{d} m(z)\right) - \max_{z\in X} g(z) \geq -\alpha + \eps \log(m(\{g-\max g \geq -\alpha\})) 
\end{equation*}
for any $\alpha \in \mathbb{R}_+$. The fact that $g$  is $L$-Lipschitz ensures that $m(\{g-\max g \geq -\alpha\})\geq m(B(x_0,\alpha/L))\geq C \left(\frac{\alpha}{L}\right)^{\delta}$ (where $x_0$ is a maximum point of $g$ and $C$ is a positive constant obtained from the condition \eqref{ahlfors}, hence only depending on $m$). Thus we have
\begin{equation*}
     \eps \log\left(\int_X e^{g(z)/\eps} \mbox{d} m(z)\right) - \max_{z\in X} g(z) \geq -\alpha + \eps \Big(\delta \log\Big(\frac{\alpha}{L}\Big)+ \log(C)\Big). 
\end{equation*}
Taking $\alpha = L \eps$ grants the result since the optimal classifiers for \eqref{upvalsoftmaxreg} are uniformly Lipschitz as $\eps \to 0$ as we have shown in Proposition \ref{hreglip}.
\end{proof}

\subsection{Numerical Results}   
The interest of the regularized problem lies in the fact that under suitable differentiability assumptions on $l_i$ one  can directly use a gradient descent algorithm to compute the solution. In this section, we present numerical solutions of the regularized problem \eqref{upvalsoftmaxreg} for different values of $\eps$. In all the examples, the space $X$ will be the one dimensional torus\footnote{conveniently identified with the unit circle or with the unit interval $[0,1]$ with periodic boundary conditions.} $X=\R/\mathbb{Z}$; the cost will be functions of the distance $d$ on the torus and we will take $c_{1}=c_{-1} = c$. The loss functions will be given by $l_i(s) = \log(1+e^{-is})$. The reference measure is $m = \sum_{k=0}^{N-1} \frac{1}{N}\delta_{\frac{k}{N}}$. Throughout this section, we take $N = 10^3$ and the measures $\mu_i$ will have a density with respect to $m$.\\
\begin{figure}\caption{Impact of $\eps$}
     \captionsetup[subfigure]{justification=Centering}
     \begin{subfigure}[t]{0.45\textwidth}
         \includegraphics[width=\textwidth]{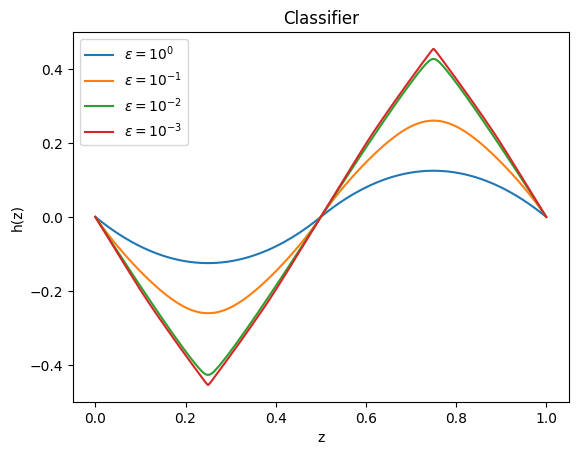}
         \caption{Graph of the classifier for varying values of $\eps$.}\label{fig:classepsilon}
     \end{subfigure}\hspace{\fill} 
     \begin{subfigure}[t]{0.45\textwidth}
         \includegraphics[width=\linewidth]{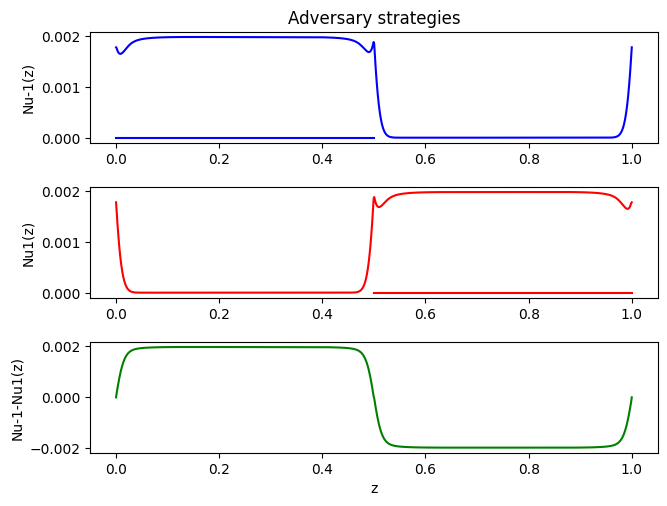}
         \caption{The adversarial strategies for $\eps= 10^{-1}$}
     \end{subfigure}
     \bigskip 
     \begin{subfigure}[t]{0.45\textwidth}
         \includegraphics[width=\linewidth]{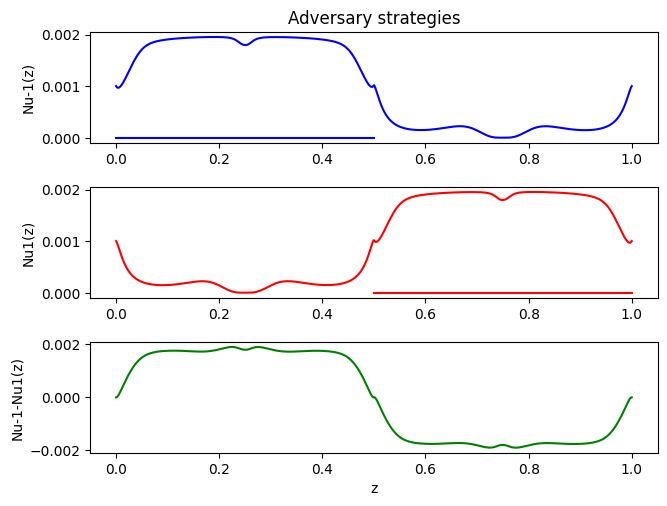}
         \caption{The adversarial strategies for $\eps= 10^{-2}$}
     \end{subfigure}\hspace{\fill} 
     \begin{subfigure}[t]{0.45\textwidth}
         \includegraphics[width=\linewidth]{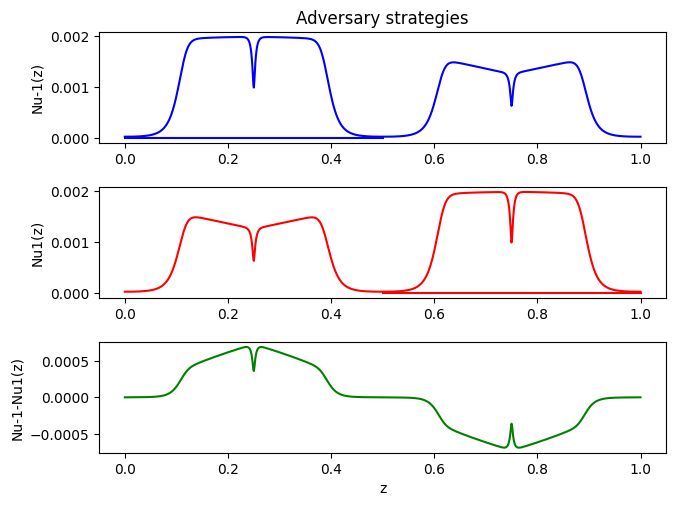}
         \caption{The adversarial strategies for $\eps= 10^{-3}$}
     \end{subfigure}
\end{figure}
\textbf{Impact of $\eps$.} In this paragraph, we look at the impact of $\eps$. We set $\mu_1$ the measure supported on $[0.5,1]$ with uniform density with respect to $m$ and $\mu_{-1}$ the measure supported over $[0,0.5]$ with uniform density with respect to $m$. The cost $c$ is taken to be equal to the distance on the torus. In Figure \ref{fig:classepsilon}, we show the classifier we obtain for different values of $\eps$. As $\eps$ decreases the classifier becomes less and less regular and more and more confident. This is a direct consequence of the fact that the larger the $\eps$ the more the adversary will be able to move mass around. Note that the singularity of the optimal adversarial attack, which involves splitting the mass, only appears for $\eps = 10^{-3}$.\\
\begin{figure}\caption{Impact of the cost}
     \captionsetup[subfigure]{justification=Centering}
     \begin{subfigure}[t]{0.45\textwidth}
         \includegraphics[width=\textwidth]{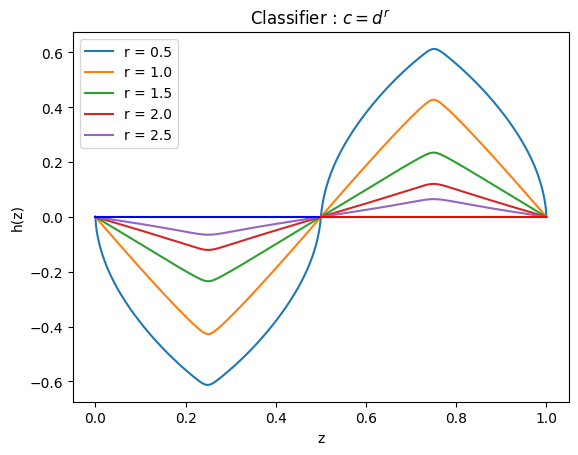}
         \caption{Graph of the classifier for varying values of the cost $c= d^r$.}\label{fig:classcostpuiss}
     \end{subfigure}\hspace{\fill} 
     \begin{subfigure}[t]{0.45\textwidth}
         \includegraphics[width=\linewidth]{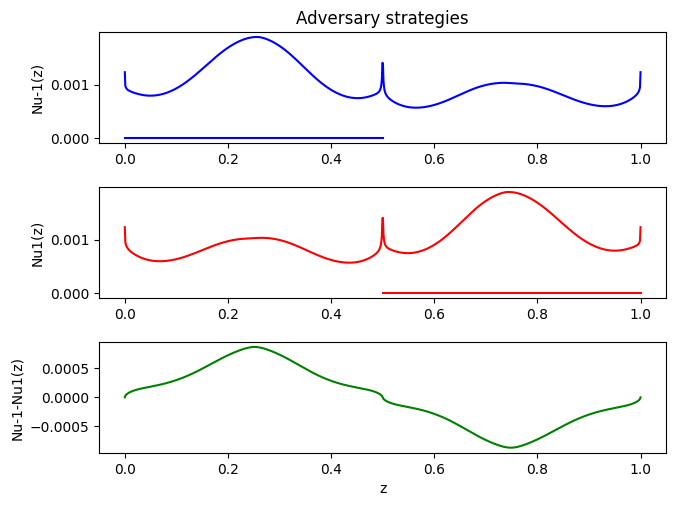}
         \caption{The adversarial strategies for $r=0.5$}
     \end{subfigure}
     \bigskip 
     \begin{subfigure}[t]{0.45\textwidth}
         \includegraphics[width=\linewidth]{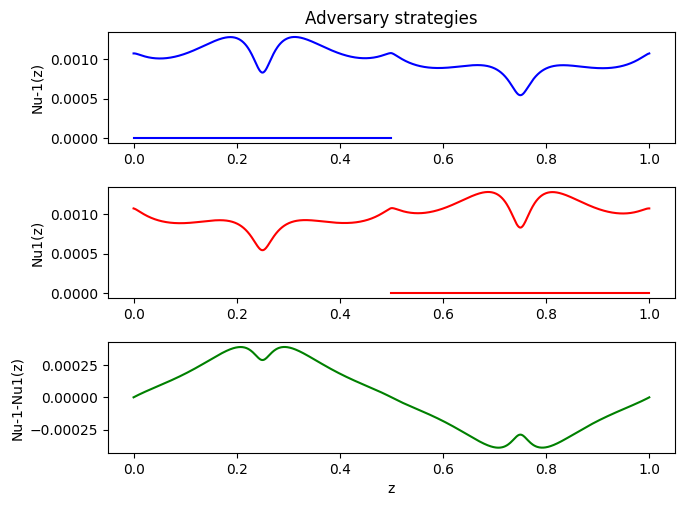}
         \caption{The adversarial strategies for $r=1$}
     \end{subfigure}\hspace{\fill} 
     \begin{subfigure}[t]{0.45\textwidth}
         \includegraphics[width=\linewidth]{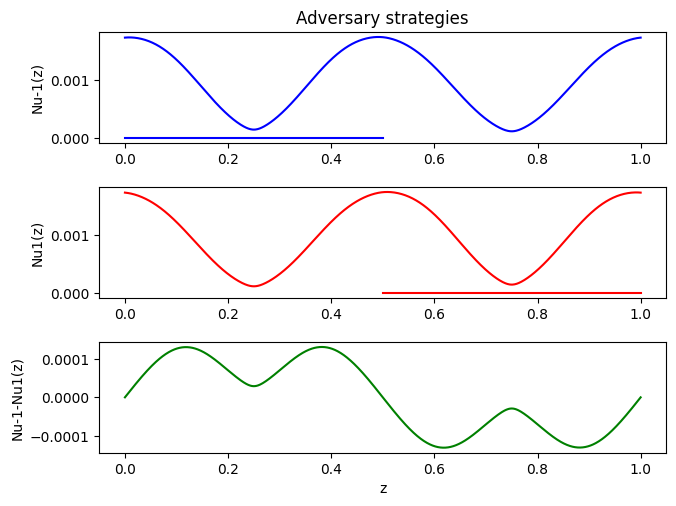}
         \caption{The adversarial strategies for $r= 1.5$}
     \end{subfigure}
\end{figure}
\textbf{Impact of the cost.} We now look at the impact of the power $r$ on the distance. The measures are the same as in the study of the impact of $\eps$. The regularity parameter $\eps$ is taken to be equal to $10^{-2}$. In Figure \ref{fig:classcostpuiss}, we display the classifer we obtain for different costs of the form $c(x,z) = d(x,z)^r$ for $r \in \{0.5,1,1.5,2,2.5\}$. Notice that as the power $r$ increases it becomes less and less costly for the adversary to move mass arround. This is coherent with the fact the absolute value of the classifier $h$ decreases with $r$. Note as well that as $r$ increases $\nu_1$ and $\nu_{-1}$ become closer and closer. It is also worth noting the cusp of the optimal attack for the square root cost.\\
\begin{figure}\caption{Impact of the initial measures}
     \captionsetup[subfigure]{justification=Centering}
     \begin{subfigure}[t]{0.45\textwidth}
         \includegraphics[width=\textwidth]{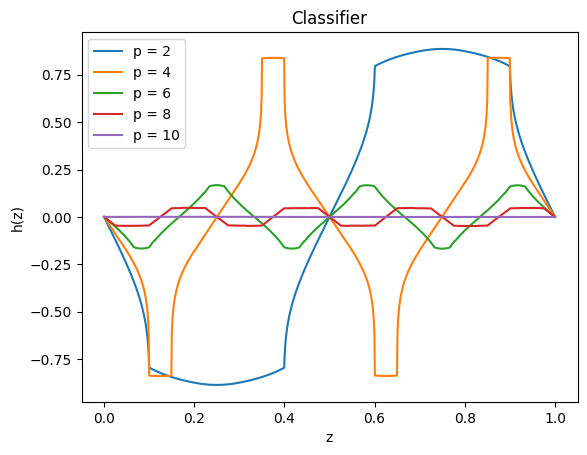}
         \caption{Graph of the classifier for varying initial measures $\mu_i$.}\label{fig:classmeasures}
     \end{subfigure}\hspace{\fill} 
     \begin{subfigure}[t]{0.45\textwidth}
         \includegraphics[width=\linewidth]{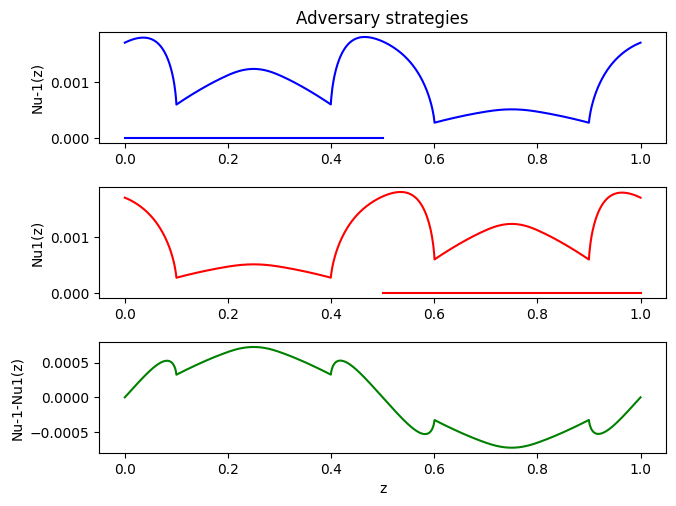}
         \caption{The adversarial strategies for $p=2$}
     \end{subfigure}
     \bigskip 
     \begin{subfigure}[t]{0.45\textwidth}
         \includegraphics[width=\linewidth]{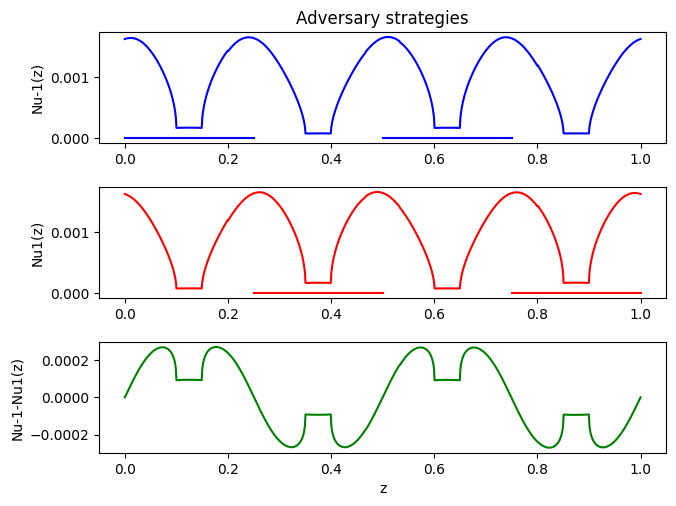}
         \caption{The adversarial strategies for $p=4$}
     \end{subfigure}\hspace{\fill} 
     \begin{subfigure}[t]{0.45\textwidth}
         \includegraphics[width=\linewidth]{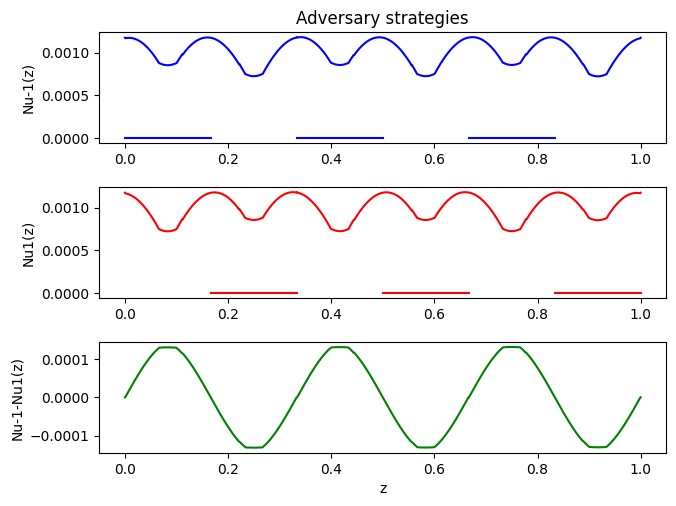}
         \caption{The adversarial strategies for $p= 6$}
     \end{subfigure}
\end{figure}
\textbf{Impact of the initial measures.} Finally, we take a discontinuous cost which prevents the adversary from moving mass too far away, namely $c(x,z) = 1_{d(x,z) > 0.1}$. The measures $\mu_i$ are taken to be constant over interleaving intervals of growing number $p\in \{2,4,6,8,10\}$. More precisely, $\mu_{-1}$ is uniform over $\cup_{k=0}^{p/2}[2k/p,(2k+1)/p]$ and $\mu_{1}$ is uniform over the complement. Figure \ref{fig:classmeasures} shows the evolution of the classifiers with the total number of intervals. When there are $10$ intervals it is possible for the adversary to match $\nu_1$ and $\nu_{-1}$ without incurring any cost. This explains that the classifier cannot distinguishs between the two measures.\\

\section{Conclusions}\label{sec-concl}

In this work, we  have considered a zero-sum game between a player designing a soft binary classifier and an adversary who may corrupt the true distributions at the expense of some transport cost. Under some mild regularity assumptions on the loss and transport cost functions, we have shown that this game has a value. More importantly, in our opinion, we established the existence of \emph{continuous} optimal classifiers by exploiting some specific  properties of  a convex minimization problem (see \eqref{upvalsimpl})  sharing some similarities with the Kantorovich dual formulation of optimal transportation. We also proposed a softmax approximation of this problem, studied its convergence, as the regularization parameter vanishes, and presented some numerical simulations based on this approximation. In these simulations, we  considered  toy one-dimensional situations  (not because of special properties of optimal transport in one dimension but because of the dimension of the softmax approximation which involves integrations with respect to a reference measure $m$ with full suport).

\smallskip

Among possible perspectives, let us mention, on the theoretical side, the extension of our analysis to adversarial multi-class classification \cite{GTKJ}. On the computational side, investigating efficient schemes that are able to treat more realistic higher-dimensional examples would be desirable.

\bigskip
  
{\bf Acknowledgments:}  G.C. acknowledges the support of the Lagrange Mathematics and Computing Research Center.

\bibliographystyle{plain}
\bibliography{bibli}

\begin{thebibliography}{10}

\bibitem{arjovsky}
Martin Arjovsky, Soumith Chintala, and L{\'e}on Bottou.
\newblock {W}asserstein generative adversarial networks.
\newblock In Doina Precup and Yee~Whye Teh, editors, {\em Proceedings of the 34th International Conference on Machine Learning}, volume~70 of {\em Proceedings of Machine Learning Research}, pages 214--223. PMLR, 06--11 Aug 2017.

\bibitem{Malick}
Wa\"iss Azizian, Franck Iutzeler, and J\'er\^ome Malick.
\newblock Regularization for {W}asserstein distributionally robust optimization.
\newblock {\em ESAIM Control Optim. Calc. Var.}, 29:Paper No. 33, 31, 2023.

\bibitem{BertsekasShreve}
Dimitri~P. Bertsekas and Steven~E. Shreve.
\newblock {\em Stochastic optimal control}, volume 139 of {\em Mathematics in Science and Engineering}.
\newblock Academic Press, Inc. [Harcourt Brace Jovanovich, Publishers], New York-London, 1978.
\newblock The discrete time case.

\bibitem{Blanchet}
Jose Blanchet, Yang Kang, and Karthyek Murthy.
\newblock Robust {W}asserstein profile inference and applications to machine learning.
\newblock {\em J. Appl. Probab.}, 56(3):830--857, 2019.

\bibitem{Braides}
Andrea Braides.
\newblock {\em {$\Gamma$}-convergence for beginners}, volume~22 of {\em Oxford Lecture Series in Mathematics and its Applications}.
\newblock Oxford University Press, Oxford, 2002.

\bibitem{cuturilightspeed}
Marco Cuturi.
\newblock Sinkhorn distances: Lightspeed computation of optimal transport.
\newblock {\em Advances in Neural Information Processing Systems}, 26, 2013.

\bibitem{cuturifast}
Marco Cuturi and Arnaud Doucet.
\newblock Fast computation of wasserstein barycenters.
\newblock In {\em International conference on machine learning}, pages 685--693. PMLR, 2014.

\bibitem{gerodim}
Simone Di~Marino and Augusto Gerolin.
\newblock An optimal transport approach for the {S}chr\"odinger bridge problem and convergence of {S}inkhorn algorithm.
\newblock {\em J. Sci. Comput.}, 85(2):Paper No. 27, 28, 2020.

\bibitem{EkelandTemam}
Ivar Ekeland and Roger T\'{e}mam.
\newblock {\em Convex analysis and variational problems}, volume~28 of {\em Classics in Applied Mathematics}.
\newblock Society for Industrial and Applied Mathematics (SIAM), Philadelphia, PA, english edition, 1999.
\newblock Translated from the French.

\bibitem{GT1}
Nicol\'as Garc\'ia~Trillos, Jakwang Kim, and Matt Jacobs.
\newblock The multimarginal optimal transport formulation of adversarial multiclass classification.
\newblock {\em J. Mach. Learn. Res.}, 24:Paper No. [45], 56, 2023.

\bibitem{GTKJ}
Nicol\'{a}s Garc\'{\i}a~Trillos, Jakwang Kim, and Matt Jacobs.
\newblock The multimarginal optimal transport formulation of adversarial multiclass classification.
\newblock {\em J. Mach. Learn. Res.}, 24:Paper No. [45], 56, 2023.

\bibitem{NIPSGoodfellow2014}
Ian Goodfellow, Jean Pouget-Abadie, Mehdi Mirza, Bing Xu, David Warde-Farley, Sherjil Ozair, Aaron Courville, and Yoshua Bengio.
\newblock Generative adversarial nets.
\newblock In Z.~Ghahramani, M.~Welling, C.~Cortes, N.~Lawrence, and K.Q. Weinberger, editors, {\em Advances in Neural Information Processing Systems}, volume~27. Curran Associates, Inc., 2014.

\bibitem{GoodfellowHarness}
Ian~J. Goodfellow, Jonathon Shlens, and Christian Szegedy.
\newblock Explaining and harnessing adversarial examples.
\newblock In Yoshua Bengio and Yann LeCun, editors, {\em 3rd International Conference on Learning Representations, {ICLR} 2015, San Diego, CA, USA, May 7-9, 2015, Conference Track Proceedings}, 2015.

\bibitem{Leonard}
Christian L\'{e}onard.
\newblock From the {S}chr\"{o}dinger problem to the {M}onge-{K}antorovich problem.
\newblock {\em J. Funct. Anal.}, 262(4):1879--1920, 2012.

\bibitem{madry}
Aleksander Madry, Aleksandar Makelov, Ludwig Schmidt, Dimitris Tsipras, and Adrian Vladu.
\newblock Towards deep learning models resistant to adversarial attacks.
\newblock In {\em International Conference on Learning Representations}, 2018.

\bibitem{Meunieretal}
Laurent Meunier, Meyer Scetbon, Rafael~B Pinot, Jamal Atif, and Yann Chevaleyre.
\newblock Mixed nash equilibria in the adversarial examples game.
\newblock In Marina Meila and Tong Zhang, editors, {\em Proceedings of the 38th International Conference on Machine Learning}, volume 139 of {\em Proceedings of Machine Learning Research}, pages 7677--7687. PMLR, 18--24 Jul 2021.

\bibitem{Kuhn}
Peyman Mohajerin~Esfahani and Daniel Kuhn.
\newblock Data-driven distributionally robust optimization using the {W}asserstein metric: performance guarantees and tractable reformulations.
\newblock {\em Math. Program.}, 171(1-2):115--166, 2018.

\bibitem{Nutz}
Marcel Nutz.
\newblock Introduction to {E}ntropic {O}ptimal {T}ransport, 2021.

\bibitem{CuturiPeyre}
Gabriel Peyr{\'e} and Marco Cuturi.
\newblock Computational optimal transport: With applications to data science.
\newblock {\em Foundations and Trends{\textregistered} in Machine Learning}, 11(5-6):355--607, 2019.

\bibitem{pinot}
Rafael Pinot, Raphael Ettedgui, Geovani Rizk, Yann Chevaleyre, and Jamal Atif.
\newblock Randomization matters how to defend against strong adversarial attacks.
\newblock In Hal~Daum\'e III and Aarti Singh, editors, {\em Proceedings of the 37th International Conference on Machine Learning}, volume 119 of {\em Proceedings of Machine Learning Research}, pages 7717--7727. PMLR, 13--18 Jul 2020.

\bibitem{Rudin}
Walter Rudin.
\newblock {\em Real and complex analysis}.
\newblock McGraw-Hill Book Co., New York, third edition, 1987.

\bibitem{Santambrogio}
Filippo Santambrogio.
\newblock {\em Optimal transport for applied mathematicians}, volume~87 of {\em Progress in Nonlinear Differential Equations and their Applications}.
\newblock Birkh\"{a}user/Springer, Cham, 2015.
\newblock Calculus of variations, PDEs, and modeling.

\bibitem{GoodfellowIntrig}
Christian Szegedy, Wojciech Zaremba, Ilya Sutskever, Joan Bruna, Dumitru Erhan, Ian~J. Goodfellow, and Rob Fergus.
\newblock Intriguing properties of neural networks.
\newblock In Yoshua Bengio and Yann LeCun, editors, {\em 2nd International Conference on Learning Representations, {ICLR} 2014, Banff, AB, Canada, April 14-16, 2014, Conference Track Proceedings}, 2014.

\bibitem{GTexist}
Nicolas~Garcia Trillos, Matt Jacobs, and Jakwang Kim.
\newblock On the existence of solutions to adversarial training in multiclass classification, 2023.

\bibitem{GT2}
Nicolas~Garcia Trillos, Matt Jacobs, Jakwang Kim, and Matthew Werenski.
\newblock An optimal transport approach for computing adversarial training lower bounds in multiclass classification, 2024.

\bibitem{Villani}
C\'{e}dric Villani.
\newblock {\em Topics in optimal transportation}, volume~58 of {\em Graduate Studies in Mathematics}.
\newblock American Mathematical Society, Providence, RI, 2003.

\end{thebibliography}

\end{document}